\documentclass[12pt]{amsart}

\usepackage{amsmath,amsfonts,amsthm,amsopn,cite,mathrsfs}
\usepackage{epsfig,verbatim}
\usepackage{subfigure}

\newcommand{\tfa}{time-frequency analysis}

\newcommand{\ft}{Fourier transform}
\newcommand{\stft}{short-time Fourier transform}

\newcommand{\tf}{time-frequency}

\newcommand{\fif}{if and only if}
\newcommand{\tfs}{time-frequency shift}

\newcommand{\modsp}{modulation space}
\newcommand{\psdo}{pseudodifferential operator}

\newtheorem{tm}{Theorem}
\newtheorem{lemma}[tm]{Lemma}
\newtheorem{prop}[tm]{Proposition}
\newtheorem{cor}[tm]{Corollary}

\newcommand{\rems}{\noindent\textsl{REMARKS:}}
\newcommand{\rem}{\noindent\textsl{REMARK:}}

 \theoremstyle{definition}
 \newtheorem{definition}{Definition}

\newcommand{\beqa}{\begin{eqnarray*}}
\newcommand{\eeqa}{\end{eqnarray*}}

\newcommand{\field}[1]{\mathbb{#1}}
\newcommand{\bR}{\field{R}}        
\newcommand{\bZ}{\field{Z}}        
\newcommand{\vf}{\varphi}

 \def\cS{\mathcal{S}}
 
 \def\cH{\mathcal{H}}

 \def\cG{\mathcal{G}}

 \def\cA{\mathcal{A}}

 \def\cO{\mathcal{O}}

\def\rd{\bR^d}

\def\rdd{{\bR^{2d}}}
\def\zdd{{\bZ^{2d}}}

\def\lrd{L^2(\rd)}

\def\mvv{M_v^1}
\def\mif{M^{\infty,1}}

\def\Lmpq{L ^{p,q}_m}

\def\Mmpq{M_m^{p,q}}

\def\intrd{\int_{\rd}}
\def\intrdd{\int_{\rdd}}

\def\<{\left<}
\def\>{\right>}

\def\inv{^{-1}}

\def\mv1{M_v^1}

\newcommand{\bba}{\mathbf{a}}

\newcommand{\bbc}{\mathbf{c}}

\newcommand{\gmu}{Gabor multiplier}
\newcommand{\sooo}{S^0_{0,0}}
\newcommand{\sdd}{\sigma (x,D)}
\newcommand{\miff}{M^{\infty, \infty}}
\begin{document}
\begin{abstract}
We investigate a new representation of general operators by means of
sums of shifted Gabor multipliers. These representations arise by studying the
matrix of an operator with respect to a Gabor frame. Each shifted
Gabor multiplier corresponds to a side-diagonal of this matrix. This
representation is especially useful for operators whose associated
matrix possesses some off-diagonal decay. In this case one can
completely characterize the symbol class of the operator by the size 
of the symbols of the Gabor multipliers. As an application we derive
approximation theorems for \psdo s in the Sj\"ostrand class.
\end{abstract}

\title[Representation and Approximation of Pseudodifferential
Operators]{Representation and Approximation of
  Pseudodifferential Operators by Sums of  Gabor 
  Multipliers}
\author{Karlheinz Gr\"ochenig}
\address{Faculty of Mathematics \\
University of Vienna \\
Nordbergstrasse 15 \\
A-1090 Vienna, Austria}
\email{karlheinz.groechenig@univie.ac.at}
\date{}
\thanks{K.~G.~was supported by the Marie-Curie Excellence Grant
  MEXT-CT 2004-517154 and in part by the FWF grant SISE S10602.}
\dedicatory{Dedicated to Paul L.\ Butzer} 

\maketitle

\section{Introduction}


One of the approaches to understand a given  operator is to decompose (=
``analyze'') the operator into simpler operators and then  to find
approximation theorems. The meaning of ``simple'' and
``approximation'' varies with one's point of view and with the application
at hand. In this investigation we take the point of view of \tfa\ and
consider an operator as simple if it aligns well with Gabor frames.
Technically, the simple operators are \emph{Gabor 
  multipliers}. Given a point $z= (x,\xi ) \in \rdd$ in the \tf\
plane, we denote the associated \tfs\  acting on a function $f$ by 
$$
\pi (z) f(t) = e^{2\pi i \xi \cdot t} f(t-x),  \qquad  x,\xi ,t \in \rd
\, .
$$
Now fix a non-zero function $g\in \lrd $ (a so-called ``window''
function) and  a lattice $\Lambda = A\zdd $ (with an invertible $2d\times
2d$-matrix). A \gmu\ with the symbol $\bba = (a_\lambda )_{\lambda \in
  \Lambda }$ is defined informally as
\begin{equation}
  \label{eq:h1}
  M_\bba ^{g,\Lambda } f = \sum _{\lambda \in \Lambda } a_\lambda
  \langle f , \pi (\lambda )g\rangle \pi (\lambda ) g \, .
\end{equation}
\gmu s have attracted attention both in engineering and in
mathematics~\cite{BP06,DT09,FHK04,FN03,wong02}, because they provide an easy and
computationally attractive model of \tf\ masking. In \eqref{eq:h1} the Gabor
coefficients $\langle f, \pi (\lambda )g\rangle $ are a measure of the \tf\
content at $\lambda = (\lambda _1,\lambda _2) \in \rdd $, i.e., the
amplitude of the frequency $\lambda _2$ at time $\lambda _1$. The
multiplication of the Gabor coefficients with the symbol or
\emph{mask} $\bba $ amounts to an enhancement or damping of certain
\tf\ regions. In this sense, the symbol $\bba$ can be compared to the
transfer function of a discrete time-invariant system usually given  by
a convolution operator $ f   \to f \ast a$.

Whereas the basic properties (boundedness, mapping properties,
Schatten class properties)  of \gmu s are  well
understood~\cite{FN03}, it is still open how \gmu s fit into the
general picture of (pseudodifferential) operators. How big is the
class of \gmu s? How well can a given operator be approximated by a
\gmu\ or by a sum of \gmu s? Which  properties influence and determine
the accuracy of an  approximation by \gmu s? 

To address these questions, we will prove results of the following
type.

(1) We show that every  operator $A$  from $\cS (\rd )$ to $\cS
'(\rd )$ can be represented as a sum of shifted \gmu s
\begin{equation}
  \label{eq:h2}
A = \sum _{\nu \in \Lambda } \pi (\nu ) M_{\bba _\nu }  ^{g,\Lambda }
\end{equation}
for suitable windows $g$ and symbols sequences $\bba _\nu $. 
This decomposition is a new representation of operators comparable to
the Schwartz kernel theorem (which we will indeed use to deduce
\eqref{eq:h2}).

(2) We  relate the properties of the symbol of a \psdo\ and the
multiplier sequences $\bba _\nu $. One of the main results is a
characterization of the generalized Sj\"ostrand classes by the size of
the multipliers $\bba _\nu $. 

(3) We estimate the error 
$$
\|A - \sum _{|\nu |\leq N} \pi (\nu )  M_{\bba _\nu } \|_{L^2 \to L^2}
$$ 
in the operator norm on $\lrd $ and on other function spaces. In this way
we develop
an approximation theory for operators. One can draw 
an analogy to the classical theorems of Jackson for the approximation
of a continuous function by trigonometric polynomials,
e.g.,~\cite{butzer}.
Instead of continuous functions we consider bounded operators, and
instead of trigonometric polynomials we approximate with finite sums
of \gmu s.   Our main result is in the spirit of Jackson's theorem:
the rate of approximation by
sums of Gabor multipliers is directly correlated to the smoothness of  the symbol.

(4) Finally we  indicate  how \gmu s are in the modeling  time-varying
environments in wireless communication.

A related investigation of \gmu s was carried out by 
D\"orfler and Torresani~\cite{DT09}. They approximate  Hilbert-Schmidt
operators by \gmu s and sums of \gmu s with different windows and
perform a window optimization. 
Though similar in spirit, the results and methods are completely
disjoint. The Hilbert-Schmidt norm of an operator is accessible to
explicit calculations, and questions of best approximation and error
estimates can be treated by using orthogonal projections. By contrast,
there is no formula for the operator norm of an operator, and we have
to use different ideas. Our results are based on the almost
diagonalization of \psdo s with respect to Gabor frames~\cite{Gro06} and the
analogy between  matrices and operators. In this analogy a \gmu\
corresponds to a  diagonal
matrix, and the representation \eqref{eq:h2} corresponds to the
decomposition of a matrix into the sum of its side-diagonals. 
The approximation in the operator norm is also motivated by the needs
of wireless communications. There the arising operators are
perturbations of the identity operator and of convolution operators
and they can never be Hilbert-Schmidt operators.


The paper is organized as follows: In Section 2 we recall a minimum of
definitions and results from \tfa , in Section~3 we introduce the main
tool for the approximation theory by Gabor multipliers, namely the
almost diagonalization of \psdo s in the Sj\"ostrand class. Section~4
presents the formal definition and boundedness properties of \gmu
s. In Section~5 we derive several versions  of the operator 
representation \eqref{eq:h2},  and in  Section~6 we investigate the
approximation of operators by \gmu s. We conclude with  a brief discussion of \gmu s in
wireless communications in Section~7. 
 
\section{Some Time-Frequency Analysis}

 Let $z= (x,\xi ) \in \rd \times \rd $ be a point in the \tf\ plane,
$t\in \rd $. The \tfs\ $\pi (z) $ is the operator 
$$
\pi (z) f(t) = e^{2\pi i \xi \cdot t} f(t-x)
$$
acting on $L^2(\rd ), \cS ' (\rd )$ and many other spaces. In the
following we use  \cite{book} as a 
general reference to \tfa .

Associated to this set of operators are a signal transform, function
spaces, and structured frames; namely the \stft , the  \modsp s, and
Gabor frames.


 Fix a  non-zero  window function $g$,   $g\in \cS $, say. 
 The 
\emph{short-time \ft } of $f$ with respect to the ``window'' $g$ is
defined as  
$$
V_gf (z)  = \langle f, \pi (z) g\rangle = \intrd f(t)
\overline{g(t-x)} e^{2\pi i \xi \cdot t} \, dt = \big( f \cdot g(. -
x)\big) \, \widehat{} \, (\xi ) \, , 
$$
whenever the duality $\langle \cdot , \cdot \rangle $, the integral,
or the \ft\ are defined.    

Further, we may associate a class of function spaces to the \tfs s
and the \stft . Again for fixed non-zero window function $g$, $g\in
\cS $, say,  we define the 
modulation spaces $\Mmpq $ for  $1\leq p,q\leq \infty $ and $m$ a 
weight function as follows: 
$$
f\in \Mmpq (\rd ) \, \Leftrightarrow \,  V_g f \in \Lmpq 
$$
 with norm $\|f\|_{\Mmpq } = \|V_gf \|_{\Lmpq } $. Here $\Lmpq $ is
 the usual mixed-norm space. These spaces are well-defined, there is
 an extensive theory about \modsp s. The reader  should
 consult~\cite[Chps.~11-13]{book} and the references therein, an
 important source for the history is~\cite{feiSTSIP}.

Another notion involving the set of  \tfs s $\pi (z)$ is the notion of  
Gabor frames. Fix a lattice $\Lambda \subseteq \rdd $, i.e.,  $\Lambda
= A\zdd $ for an invertible  $2d\times 2d$  matrix $A$. 
The collection of \tfs s 
$\cG (g,\Lambda ) = \{ \pi (\lambda )g: \lambda \in \Lambda \}$ for
some non-zero $g\in \lrd $ is called a Gabor system. The set $\cG
(g,\Lambda )$ is  a Gabor frame, if there exist constants
$A,B>0$ such that 
\begin{equation}
  \label{eq:1}
A \|f\|_2^2 \leq \sum _{\lambda \in \Lambda } | \langle f, \pi
(\lambda ) g\rangle |^2 \leq B \|f\|_2^2  \qquad \forall f\in L^2(\rd
)   \, .
\end{equation}
 If $A=B=1$, $\cG (g,\Lambda )$ is called a  Parseval frame and
 \eqref{eq:1} implies the expansion 
 \begin{equation}
   \label{eq:heu}
f = \sum _{\lambda \in \Lambda } \langle f,\pi (\lambda ) g\rangle \pi
(\lambda ) g \qquad \forall f\in \lrd    
 \end{equation}
with unconditional convergence in $\lrd $. In contrast to an
orthonormal expansion the coefficients in this expansion need not be
unique, the choice $\langle f, \pi (\lambda ) g\rangle $ is a
distinguished, explicit, and convenient choice. The existence problem
for Gabor frames is almost completely settled~\cite{bekka04}. In what
follows, we take the existence of Gabor frames and Parseval Gabor
frames for granted. 
 
\vspace{3 mm}

\textbf{Characterizations of \modsp s by Gabor frames.}
In the following $v$ always denotes a submultiplicative, even weight
on $\rdd $,
i.e., $v(x+y) \leq v(x) v(y)$. A weight function $m$ on $\rdd $  is called
$v$-moderate, if $m(x+y) \leq C v(x) m(y)$ for all $x,y\in \rdd $. 

  If  $\cG
(g,\Lambda )$ is a Gabor frame and $g\in \mvv (\rd )$, then magnitude
of the frame
coefficients $\langle f, \pi (\lambda )g\rangle , \lambda \in \Lambda
,$ characterizes the membership of a function in a particular
modulation space~\cite{fg97jfa,book}.
 More precisely,  for
every $p\in [1,\infty ]$ and $v$-moderate weight $m$, 
$$
f \in M^{p,p}_m(\rd )  \Leftrightarrow \Big( \sum _{\lambda \in
  \Lambda } |\langle f, \pi (\lambda )g\rangle |^p m(\lambda )^p
\Big)^{1/p} < \infty \, ,
$$
and the latter expression is an equivalent norm on $M^{p,p}_m$. An
analogous characterization holds for rectangular lattices and the
mixed \modsp s $\Mmpq$~\cite{fg97jfa}.  For  $p\neq q$ 
and a  general lattice $\Lambda $, we have $f\in \Mmpq $, \fif\ 
  $\sum _{\lambda \in \Lambda } |\langle f , \pi (\lambda )g\rangle |
  \chi _{\lambda +Q} \in \Lmpq (\rdd )$\cite{fg89jfa}.

In particular, we have the following simple characterization of the
Schwartz class by means of Gabor frames. Assume that $g\in \cS (\rd )$
and that $\cG (g, \Lambda )$ is a frame for $\lrd $. A function $f \in
\lrd $ belongs to $\cS (\rd )$, \fif\  $\langle f, \pi (\lambda
)g\rangle , \lambda \in \Lambda $ decays rapidly, i.e.,
\begin{equation}
  \label{eq:cd6}
|\langle f, \pi (\lambda )g\rangle | = \cO \Big( (1+|\lambda
|)^{-s}\Big) \qquad \text{ for all } s\geq 0 \, . 
\end{equation}

\vspace{3 mm}

\textbf{Pseudodifferential Operators.} Given a function (or distribution) $\sigma $ on $\rdd $, the
corresponding 
pseudodifferential operator with symbol $\sigma $ is defined
informally by the integral
\begin{equation}
  \label{eq:2}
  \sdd f(x) = \int _{\rd } \sigma (x,\xi  ) \hat{f} (\xi ) e^{2\pi i
  x\cdot \xi } \, d\xi \, .
\end{equation} 
Again the definition of \psdo s does not necessarily require that the
integral is defined, see~\cite{Gro06ped,hormander3} for a rigorous definition. 

\vspace{3 mm}

\textbf{Symbol Classes.} We consider certain \modsp s as symbol
classes. They are generalizations of the Sj\"ostrand
class~\cite{Sjo94} and may be understood as non-smooth extensions of
the H\"ormander class $\sooo$. Fix a non-zero test function $\Phi $ on
$\rdd $, e.g., the Gaussian $\Phi (z) = e^{-\pi z\cdot z}$ and a
submultiplicative weight $v$ on $\rdd $. Then a symbol $\sigma $
belongs to $\mif _v (\rdd )$, if 
\begin{equation}
  \label{eq:3}
  \intrdd \sup _{z\in \rdd }|V_\Phi \sigma
  (z,\zeta )|\, v(\zeta )\,  d\zeta  = \|\sigma \|_{\mif _{v}} < \infty \, .
\end{equation}
Likewise, a symbol $\sigma $ belongs to $M^{\infty , \infty }_v (\rdd
)$, if 
\begin{equation}
  \label{eq:4}
  \sup _{\zeta \in   \rdd } \sup _{z\in \rdd }|V_\Phi \sigma
  (z,\zeta )| v(\zeta ) = \|\sigma \|_{M^{\infty , \infty }_v} < \infty \, .
\end{equation}

The generalized Sj\"ostrand classes are thus certain \modsp s on $\rdd
$. Further generalizations use the class of solid convolution
algebras on $\zdd $ to parametrize the extensions of $\sooo
$~\cite{GR08}.

\section{Almost Diagonalization}

We will  study the matrix of a \psdo\ $\sdd $ with respect to a Gabor
frame. 

Precisely, let $\cG(g,\Lambda )$ be a  Gabor system  and $\sigma
$ a symbol with corresponding \psdo\ $\sdd $. Then we 
consider the matrix $M =
M(\sigma )$ with entries
\begin{equation}
  \label{eq:cd8}
M(\sigma ) _{\lambda , \mu } = \langle \sdd (\pi (\mu ) g), \pi
(\lambda ) g \rangle  \qquad \lambda ,\mu \in \Lambda \, .  
\end{equation}

The main results of \cite{Gro06,GR08} establish a precise link between
\psdo s in  generalized Sj\"ostrand classes and the associated
matrix. In the following $j$ denotes the rotation $j(z_1,z_2)= (z_2,
-z_1)$. (It is needed to use the standard modulation operators instead
of some symplectic modulations.)

\begin{tm}[\cite{Gro06}] \label{sparsity}
Assume that $g\in M^1 _{v }  (\rd )$ and that $\cG
(g,\Lambda )$ is a  Gabor
frame. Then  $\sigma \in \mif _{v\circ j\inv } $ \fif\ 
  there exists  $h \in \ell ^1 _v (\Lambda)$ such that 
  \begin{equation}
    \label{eq:do2}
|M(\sigma ) _{\lambda,\mu}| = |\langle \sdd  \pi (\mu )g,
\pi (\lambda) g \rangle| \leq h(\lambda -\mu  )    
  \end{equation}
for all $\lambda , \mu \in \Lambda $. Furthermore, $\inf \|h\|_{\ell
  ^1_v}$ with the infimum taken over all $h$ satisfying \eqref{eq:do2}
is an equivalent norm on $\mif _{v\circ j\inv} $.  
\end{tm}

\rem\ If $\cG (g, \Lambda )$ is a merely a Gabor system with $g\in M^1_{v\circ
  j\inv}(\rd )$ (but not a frame), then the almost diagonalization \eqref{eq:do2} still
holds. However, in general, the converse is not true in general,
because the matrix $M(\sigma )$ does not fully describe the operator
$\sdd $.  
 
A similar theorem holds for the class $M^{\infty , \infty }_v$. 
\begin{tm}[\cite{GR08}] \label{sparsity2}
Assume that $g\in M^1 _{v }  (\rd )$ and that $\cG
(g,\Lambda )$ is a  Gabor
frame. Then  $\sigma \in \miff _{v\circ j\inv } $ \fif\ 
  there exists  
  \begin{equation}
    \label{eq:do1}
|M(\sigma ) _{\lambda,\mu}| = |\langle \sdd  \pi (\mu )g,
\pi (\lambda) g \rangle| \leq C v(\lambda -\mu  )\inv 
  \end{equation}
for all $\lambda , \mu \in \Lambda $. 
\end{tm}
Thus in the case of $\miff _v$ one obtains genuine off-diagonal decay
of the associated matrix $M(\sigma )$. The standard weights to
describe the decay condition are either polynomial weights $v(z)  =
\langle z\rangle ^s = (1+|z|^2)^{s/2}$ or (sub)exponential weights $v(z)
= e^{a|z|^b}$ for $a>0$ and $0\leq b \leq 1$.

Other types of decay conditions are studied in ~\cite{GR08}.

As a special case we mention the  H\"ormander Class $\sooo $
consisting of all $C^\infty $-functions on $\rdd $ with bounded derivatives. This
symbol class is related to \modsp s as follows~\cite{GR08,Toft02a}:
$\sooo = \bigcap _{s>0} \miff _{\langle\zeta \rangle ^s}$. As a
consequence of Theorem~\ref{sparsity2} we obtain a characterization of
$\sooo $. 

 \begin{cor}
   $g\in \cS , g\neq 0$, $\cG (g, \Lambda )$ Gabor frame. Then 
$\sigma \in \sooo $ \fif\ for every $s\geq 0$ there is a constant $C_s$ such that 
$$
 |\langle \sdd  \pi (\mu )g, \pi (\lambda) g \rangle| \leq C_s\, (1+|\lambda
 -\mu  |)^{-s}  \qquad \text{ for all } \lambda ,\mu \in \Lambda \, .
$$ 
 \end{cor}

For completeness we mention that  operators with symbols in one of the
generalized Sj\"ostrand classes
$\mif _v$ and $\miff _v$ are closed under composition and that they
are bounded on many \modsp s~\cite{Gro06,GR08,Sjo95}.

\section{Gabor Multipliers}

Next we introduce a special class of operators, so-called Gabor
multipliers. These operators are particularly simple and, in some
sense, correspond to diagonal matrices.

\begin{definition}\label{gabmult}
Let $\cG (g,\Lambda )$ be a Gabor system. 
  Given ``symbol sequence''  $\mathbf{a}= (a_\lambda )_{\lambda \in
    \Lambda}$, the Gabor multiplier   $M_\mathbf{a}$ is defined to be
  the operator 
$$
M_\mathbf{a}f = \sum _{\mu \in \Lambda } a_\mu \langle f, \pi
(\mu )g \rangle \pi (\mu )g \, .
$$
\end{definition}
Clearly, the definition also depends on the Gabor system $\cG (g,
\Lambda )$. To indicate the dependence of the Gabor multiplier of all
parameters, we would have to write $M_\bba ^{g,\Lambda }$. Since our
results holds generically for all Gabor systems in a certain class, we
prefer to keep the notation simple and 
 omit the reference to $\cG (g,\Lambda )$.  Likewise one could use two
 windows $g$ and $\gamma $ and consider \gmu s of the form 
$M_\mathbf{a}f = \sum _{\mu \in \Lambda } a_\mu \langle f, \pi
(\mu )g \rangle \pi (\mu )\gamma $. Such generalizations cause  only notational
complications, but do not change any of the results or arguments. 

Gabor multipliers have been studied in detail by Feichtinger and
Nowak~\cite{FN03}. Gabor multipliers can also be interpreted as \tf\  localization
operators with distributional symbols, and boundedness results  follow
from the theory of localization operators~\cite{CG03}. 

As a typical boundedness result we state the following one. 

\begin{lemma}
  Assume that $g\in \mvv (\rd )$. If $\bba \in \ell ^\infty (\Lambda
  )$, then $M_\bba $ is bounded on every \modsp\ $\Mmpq $ for $1\leq
  p,q\leq \infty $ and every $v$-moderate weight $m$ with a uniform bound for
  the operator norm 
  \begin{equation}
    \label{eq:h3}
    \|M_\bba \|_{\Mmpq \to \Mmpq } \leq C_\Lambda \|g \|_{\mvv }^2 \,
    \|\bba \|_{\infty } \, .
  \end{equation}
\end{lemma}
\begin{proof}
  The decisive property for the proof is the boundedness of the
  coefficient operator $f \to \Big( \langle f, \pi (\lambda )g\rangle
  \Big)_{\lambda \in \Lambda }$ from $\Mmpq (\rd ) $ to the sequence space
  $\ell _m^{p,q}(\Lambda )$ and the synthesis operator $\bbc \to \sum _{\lambda
    \in \Lambda } c_\lambda \pi (\lambda
  )g$. See~\cite[Thms.~12.2.1-4]{book} or \cite{fg97jfa} (for
  rectangular lattices).    
Using these estimates, we find for $p=q$ that 
\begin{eqnarray}
\|M_\mathbf{a}f\|_{M^{p,p}_m}  &=& \|\sum _{\mu \in \Lambda } a_\mu \langle f, \pi
(\mu )g \rangle \pi (\mu )g\|_{M^{p,p}_m} \notag \\
& \leq &  C \|g\|_{\mvv } \big( \sum _{\mu \in \Lambda }   | a_\mu \langle f, \pi
(\mu )g \rangle |^p m(\mu )^p \big)^{1/p} \notag \\
&\leq & C \|g\|_{\mvv } \, \|\bba \|_\infty \big(\sum |\langle f, \pi
(\mu )g \rangle |^p m(\mu )^p \big)^{1/p}\label{eq:do3} \\
&\leq & (C \|g\|_{\mvv })^2 \, \|\bba \|_{\infty } \|f\|_{M^{p,p}_m}  \, . \notag
\end{eqnarray}
For $p\neq q$ and non-separable lattices the proof is identical, once
the sequence space is defined correctly. See~\cite{fg89jfa}.
\end{proof}


\section{Representations of Operators with Gabor Multipliers}

After these preparations let us now explain 
   how and why Gabor multipliers arise in the theory of \psdo s. 
Assume that $\cG (g, \Lambda )$ is a Parseval frame, so that
\eqref{eq:heu} holds. We
   expand both $f$ and $Af= \sdd f $ with 
   respect to the frame $\cG (g,\Lambda )$. Then 
$ f= \sum _{\mu \in \Lambda  }  \langle f, \pi (\mu )
g\rangle \pi (\mu  ) g$ and  $\sdd f=\sum _{\lambda \in \Lambda }
\langle Af, \pi (\lambda )g\rangle \pi 
  (\lambda )g $. Substituting the expansion of $f$ into the
  coefficients $\langle Af, \pi (\lambda )g\rangle$, we obtain that   
\begin{eqnarray}
  Af &=& \sum _{\lambda \in \Lambda } \langle Af, \pi (\lambda )g\rangle \pi
  (\lambda )g \notag \\  
&=&  \sum _{\lambda \in \Lambda } \Big(  \sum _{\mu \in \Lambda }
\langle A    \pi (\mu )g, \pi (\lambda
)g\rangle \langle f, \pi   (\mu )g \rangle \Big) \, \pi (\lambda )g
\qquad \qquad \lambda  = \mu +\nu   \notag
\\   
&=&  \sum _{\nu \in \Lambda } \Big(  \sum _{\mu \in \Lambda } \langle
A \pi (\mu )g, \pi (\mu +\nu)g\rangle \langle f, \pi   (\mu )g \rangle
\Big) \, \pi (\mu  +\nu )g \notag
\\   
 &=& \sum _{\nu \in \Lambda } \pi (\nu ) \, \Big(  \sum _{\mu \in
   \Lambda } \langle A \pi (\mu )g, \pi (\mu +\nu)g\rangle \, e^{2\pi
   i \nu _1\cdot \mu _2}\,  \langle f, \pi   (\mu )g \rangle \, \pi
 (\mu   )g  \Big) \label{eq:5}
\end{eqnarray}
In the transition to the last line we have used the commutation rule   
for \tfs s  $\pi (\mu +\nu )
= e^{2\pi   i \nu _1\cdot \mu _2} \pi (\nu ) \pi (\mu )$. 

 This calculation is known, of course, to every student of linear
algebra; it shows how to express a linear operator by its matrix with
respect to a ``basis'', which in our case we take to be a Gabor frame.

 Now define the sequence $\bba _\nu $ by  
 \begin{equation}
   \label{eq:4a}
\bba _\nu (\mu ) =   \langle A \pi (\mu )g, \pi (\mu +\nu)g\rangle \, e^{2\pi
   i \nu _1\cdot \mu _2}\, , 
 \end{equation}
then the expression in~\eqref{eq:5} in parenthesis is just a Gabor multiplier with
symbol $\bba _\nu $, and we may rewrite  \eqref{eq:5}
as  
\begin{equation}
  \label{eq:6}
\sdd = \sum _{\nu \in \Lambda } \pi (\nu ) M_{\mathbf{a}_\nu } \, .  
\end{equation}
At least informally, the above identity shows that every \psdo\ is a
sum of shifted Gabor multipliers. In this representation each \gmu\
corresponds to a side-diagonal of the associated matrix $M(\sigma )$.

To make this argument precise, we
offer several  versions. We first show that every reasonable operator
can be represented as a sum of shifted Gabor multipliers. 

\begin{prop}
Let $\cG (g,\Lambda )$ be a Parseval frame and $g \in \cS (\rd
)\setminus \{0\}$. 

(i) If the sequences $\bba _\nu $ satisfy the growth conditions 
$$
|\bba _\nu (\mu )| \leq C  (1+|\mu | + |\nu
|)^N \qquad \text{for } \lambda ,\mu \in \Lambda \, $$
then $\sum _{\nu \in \Lambda } \pi (\nu ) M_{\mathbf{a}_\nu }$ defines
a continuous operator from $\cS (\rd ) $ to $\cS ' (\rd )$ and the
series converges in the weak operator topology. 

(ii) Conversely, 
assume that  $A$ is a continuous operator  from $\cS (\rd )$ to $\cS
  '(\rd )$ (with the weak$^*$-topology). Then
  $A$ possesses a  representation
  \eqref{eq:6} as a sum of shifted Gabor multipliers with symbols
  $\bba _\nu $ satisfying the growth estimate
$$
|\bba _\nu (\mu )| \leq C (1+|\mu | + |\nu |)^N
$$
for constants $C, N\geq 0$ depending only on $g$ and $A$. 
\end{prop}

\begin{proof}
(i)  Fix  $\nu \in \Lambda $  and consider the Gabor multiplier
$M_{\bba _\nu }f = $  \linebreak $= \sum _{\mu \in \Lambda } \bba _\nu (\mu ) \langle
f, \pi (\mu )g\rangle \pi (\mu )g$.  Since by hypothesis $\bba _\nu $ grows
polynomially and $\langle f, \pi (\mu )g\rangle$ decays rapidly by
\eqref{eq:cd6}, the coefficients $ \bba _\nu (\mu ) \langle
f, \pi (\mu )g\rangle$ also decay rapidly, whence the sum converges in
$\cS (\rd )$ and $M_{\bba _\nu }f \in \cS (\rd )$. Consequently the
partial sums $\sum _{|\nu |\leq M} M_{\bba _\nu }f $ are in $\cS (\rd
)$, and we only need to verify that the partial sums converge in the
weak$^*$ sense to an element in $\cS ' (\rd )$. The convergence
follows from 
\begin{eqnarray*}
\big|  \langle \sum _{|\nu | \leq M} M_{\bba _\nu } f, h \rangle \big| & \leq
& \sum _{|\nu | \leq M} |\bba _\nu (\mu )|\,  |\langle f, \pi (\mu
)g\rangle | \,  |\langle  \pi (\mu )g, h\rangle | \\
& \leq & C_\ell  \sum    _{|\nu | \leq M} (1+|\mu | + |\nu |)^N
(1+|\mu |)^{-\ell }
\end{eqnarray*}
for arbitrary $\ell \geq 0$. Hence the sum converges, or equivalently,
the infinite sum $\sum _{\nu \in \Lambda } M_{\bba _\nu }f$ converges
in the weak-$^*$ sense.

(ii) Let $A$ be continuous from $\cS (\rd ) $ to $\cS '(\rd ) $.   As suggested by informal derivation~\eqref{eq:5}, we choose the
symbols $\bba _\nu $ to be $\bba _\nu (\mu ) =   \langle A \pi (\mu
)g, \pi (\mu +\nu)g\rangle \, e^{2\pi    i \nu _1\cdot \mu _2}$.  

  By the Schwartz kernel theorem there exists a tempered distribution
  $k\in \cS ' (\rdd )$, such that $\langle Af, h \rangle _{\rd }  =
  \langle k , h \otimes \bar{f}\rangle _{\rdd }$ for $f, h \in \cS
  (\rd )$. Since the STFT of a tempered distribution grows at most
  polynomially, there exists an $N\geq 0$, such that 
$$
|\langle k , \pi (\lambda )g \otimes \overline{\pi (\mu )g}\rangle |=
|V_{g \otimes \bar{g}} k (\lambda _1,\mu _1; \lambda _2, -\mu _2)|  
\leq C (1+|\lambda | + |\mu |)^N \qquad \forall \lambda , \mu \in
\Lambda \, .
$$
We deduce the claimed growth estimate for the symbols as follows:
\begin{eqnarray*}
  |\bba _\nu (\mu )| &=& |\langle A \pi (\mu )g, \pi (\mu +\nu
  )g \rangle _{\rd }| \\
&=& |\langle k, \pi (\mu +\nu )g \otimes \overline{ \pi (\mu
  )g}\rangle _{\rdd } | \\ 
& \leq & C (1+|\mu | + |\mu +\nu |)^N \leq C' (1+|\mu | + |\nu
|)^N \, .
\end{eqnarray*}

By Step (i) the sum $ \sum _{\nu \in \Lambda } \pi (\nu ) M_{\bba _\nu
}$ is a continuous operator from $\cS (\rd )$ to $\cS ' (\rd )$, and
by the choice of the multipliers $\bba _\nu $ this operator must
coincide with the given operator $A$. 
\end{proof}

Our next version of \eqref{eq:6} deals with the representation of
\psdo s. 

\begin{lemma} \label{l:4}
Assume that $\cG (g, \Lambda )$ is a
  Parseval frame and $g\in M^1_v(\rd )$. 

(i) If   the sequence of symbols $\bba _\nu $ satisfies the condition
$$
\sum _{\nu \in \Lambda } \|\bba _\nu \|_\infty v(\nu ) < \infty \, ,
$$
then   the sum of shifted
Gabor multipliers \eqref{eq:6} converges in the operator norm on
$M^{p,q}_m $ for every $1\leq p,q\leq \infty $  and every $v$-moderate weight $m$. 

(ii) If  $\sigma \in \mif _{v\circ j\inv} (\rdd )$, then the series \eqref{eq:6}
converges  in the operator norm of $\Mmpq $ to $\sdd $, and $\sum
_{\nu \in \Lambda } \|\bba _\nu \|_\infty v(\nu ) \leq C \|\sigma
\|_{\mif _{v\circ j\inv}}$.
 
\end{lemma}

\begin{proof}
  (i) By \eqref{eq:h3} the operator norm of $M_{\bba _\nu }$ on
  $\Mmpq$ is bounded by
  $ C \|\bba _\nu \|_\infty $. Further,  the operator norm of a \tfs\
  $\pi (\nu )$ on $\Mmpq $ is bounded by $v(\nu )$. Thus,  summing
  over $\nu $, we obtain  
\begin{eqnarray*}
\|\sdd \|_{M^{p,q}_m \to M^{p,q}_m } &\leq & \sum _{\nu \in \Lambda } \|\pi
(\nu ) M_{\bba _\nu } \|_{M^{p,q}_m \to M^{p,q}_m } \\
&\leq  &    \sum _{\nu \in \Lambda } v(\nu ) \|M_{\bba _\nu } \|_{M^{p,q}_m \to
  M^{p,q}_m } \\
&\leq &   \sum _{\nu \in \Lambda } \|\bba _\nu  \|_{\infty } \, v(\nu ) <
\infty \, . 
\end{eqnarray*}
Thus the series of operators \eqref{eq:6} converges absolutely in the operator
norm on $\Mmpq $.

(ii) First note that by Theorem~\ref{sparsity}
 there exists a $h\in \ell ^1(\Lambda )$ such that $\|h\|_{\ell ^1_v}
 \leq C \|\sigma \|_{\mif _{v\circ j\inv}}$ and 
$$
|\mathbf{a}_\nu (\mu )| =  |\langle
\sdd  \pi (\mu )g, \pi (\mu +\nu)g \rangle \, e^{2\pi i \nu _1\cdot
  \mu _2}| \leq h(\nu )
$$ 
By Step (i),  \eqref{eq:6} converges absolutely in the operator
norm on $\Mmpq $,  and by the choice of symbols the limit coincides
with $\sdd $. The norm estimate follows from $ \sum _{ \nu \in \Lambda
  } \|\bba _\nu \|_\infty v(\nu ) \leq \|h\|_{\ell ^1_v}
 \leq C \|\sigma \|_{\mif _{v\circ j\inv}}$.    
\end{proof}
 
\rem\ By adjusting the conditions in the above proof, many more
versions of \eqref{eq:6} can be shown to be meaningful.


 By combining the almost diagonalization of Theorems~\ref{sparsity}
 and~\ref{sparsity2}, we now  characterize operators with symbols in  the
 generalized Sj\"ostrand class by means of sums of Gabor
 multipliers. 

 \begin{prop}\label{triv}
Fix a submultiplicative weight $v$ and a non-zero $g\in M^1_v(\rdd )$ such that
 $\cG (g,\Lambda )$ is a Parseval frame. 

(i)   A symbol $\sigma $
belongs to $\mif _{v\circ j\inv }(\rdd )$, if and only if there exist
sequences $\bba _\nu \in \ell ^\infty (\Lambda )$, such that 
\begin{equation}
  \label{eq:8}
\sdd = \sum _{\nu \in \Lambda } \pi (\nu ) M_{\mathbf{a}_\nu } \, 
\end{equation}
and 
\begin{equation}
  \label{eq:7}
\sum _{\nu\in \Lambda } \|\bba _\nu \|_\infty v(\nu ) <\infty \, . 
\end{equation}
Furthermore, the sequence of multipliers $\bba _\nu $ can be chosen
such that $\sum _{\nu \in \Lambda } \|\bba _\nu \|_\infty v(\nu ) \leq
C \| \sigma \| _{\mif _v} \leq C' \sum _{\nu \in \Lambda } \|\bba _\nu
\|_\infty v(\nu ) $. 

(ii) Assume in addition that $v\inv \ast v\inv \leq Cv\inv $ ($v$ is
subconvolutive). Then   $\sigma \in M^{\infty , \infty } _{v\circ j\inv } (\rdd )$, \fif\
$\sigma(x,D)$ possesses a representation \eqref{eq:8} with 
$$
\|\bba _\nu \|_\infty \leq C v(\nu )\inv \, .
$$
Again the multipliers $\bba _\nu $ can be chosen, such that $\sup
_{\nu \in \Lambda } \|\bba _\nu \|_\infty v(\nu )$ is an equivalent
norm on $M^{\infty , \infty }_{v\circ j\inv}$. 
    \end{prop}

    \begin{proof}
The sufficiency was shown in Lemma~\ref{l:4}. If $\sigma \in \mif _{v\circ j\inv}$, then 
both \eqref{eq:7} and \eqref{eq:8} are  satisfied.

Assume conversely that an operator $A: \cS (\rd ) \to \cS '(\rd )$ is
given as a sum of Gabor multipliers   \eqref{eq:8} with symbols
satisfying  \eqref{eq:7}. Then by the Schwartz kernel theorem  $A$ possesses a
symbol $\sigma \in \cS ' (\rdd )$ and $A= \sdd $. 

To show that $\sigma $ is in $\mif _v(\rdd )$,  we  estimate the size
of entries  $M(\sigma
)_{ \lambda \mu }$ with respect to a reference frame $\cG (\vf ,
\Lambda )$  and then apply the characterization of
Theorem~\ref{sparsity}.

Fix a Gabor frame
$\cG (\vf , \Lambda )$ with $\vf \in \cS (\rd )\cap \mvv (\rd )$. Then the matrix
entries $M(\sigma )_{\lambda \mu } = \langle \sdd (\pi (\mu ) \vf ), \pi (\lambda
) \vf \rangle $  are well-defined.
 Since $\sum _{\nu \in \Lambda } M_{\bba _\nu } $ converges weakly, we
 may interchange the brackets $\langle \cdot , 
\cdot \rangle$ with  the summation over $\nu $ and obtain 
\begin{eqnarray}
  M(\sigma )_{\lambda \mu } &=&  \langle \sdd \, \pi (\mu ) \vf , \pi (\lambda
) \vf \rangle \notag \\
&=& \sum _{\nu \in \Lambda } \langle \pi (\nu ) M_{\bba _\nu } \pi
(\mu )\vf , \pi (\lambda ) \vf \rangle \label{hm2} \\
&=& \sum _{\nu \in \Lambda }  \sum _{\kappa \in \Lambda } \bba
_\nu (\kappa )  \langle \pi (\mu ) \vf , \pi (\kappa )g\rangle\,  \langle  \pi
(\nu ) \pi (\kappa )g , \pi (\lambda ) \vf \rangle   \, . \notag 
\end{eqnarray}
 Since $g\in \mvv $ and $\vf \in \cS (\rd )\cap \mvv (\rd )$, the sequence 
$$
h(\lambda ) = |\langle g, \pi (\lambda )\vf \rangle |
$$
belongs to $\ell ^1_v(\Lambda )$ and $\|h\|_{\ell ^1_v} \leq C \|g
\|_{\mvv }$ (with a constant depending only on $\vf $ and $\Lambda
$)~\cite[Prop.~12.1.11]{book}. Set $\alpha (\nu ) = \|\bba _\nu
\|_\infty $, $h^*(\lambda )= h(-\lambda )$,  take absolute
values in \eqref{hm2},  and use the decay of $h$. Then  we obtain
that 
\begin{eqnarray*}
  |M(\sigma )_{\lambda \mu }|& \leq & \sum _{\nu \in \Lambda } \sum
  _{\kappa \in \Lambda } \|\bba _\nu \|_\infty \, h(\mu - \kappa ) \,
  h(\lambda - \nu -\kappa ) \\
&=& \sum _{\nu \in \Lambda } \alpha (\nu )  \, \sum _{\kappa
  \in \Lambda } h^*(\kappa ) h(\lambda -\nu - (\kappa +\mu )) \\
&=& \sum _{\nu \in \Lambda }  \alpha (\nu ) \, (h^* \ast
h)(\lambda - \mu -\nu ) \\
&=& (\alpha \ast h^* \ast h) (\lambda -\mu )  \, .
\end{eqnarray*}
Since $\alpha \in \ell ^1_v$ by assumption~\eqref{eq:8} and $h\in \ell
^1_v$ because $g\in \mvv $, the matrix $M(\sigma )$ is dominated
entrywise by the sequence $\alpha \ast h^* \ast h \in \ell ^1_v$. 

 Theorem~\ref{sparsity} applies and we conclude
that $\sigma \in M^{\infty , 1} _{v\circ j\inv} (\rdd )$, as claimed. Furthermore, 
\begin{eqnarray*}
  \|\sigma \|_{\mif _v} &\leq & C \|\alpha \ast h^* \ast h \|_{ \ell
    ^1_v} \\
&\leq & C \|\alpha \|_{\ell ^1_v} \|h\|_{\ell _v^1}^2 \\
&\leq & C ' \|\alpha \|_{\ell ^1_v} \|g\|_{\mvv } ^2 \, .
\end{eqnarray*}

The proof for $\sigma \in M^{\infty , \infty } _{v\circ j\inv}$ is
similar. The only modification occurs in the last part, where we have
to use the subconvolutivity of $v$ and  the convolution relation $\ell
^\infty _v \ast \ell ^\infty _v \subseteq \ell ^\infty _v$ 
\end{proof}

\section{Approximation Theorems}

Since every operator $A$ can be represented as an infinite sum of shifted
Gabor multipliers, it is natural to truncate the infinite series
\eqref{eq:8} and approximate $A$ by a finite sum of Gabor
multipliers. 
Compared to a general operator, a finite sum of Gabor multipliers is
easy to understand and easy to treat computationally. 

For the formulation of the approximation theorems we introduce the error
\begin{equation}
  \label{eq:9}
  E_N(\sigma ) := \| \sdd - \sum _{|\nu
  |\leq N} \pi (\nu ) M_{\mathbf{a}_\nu }\| _{M^{p,q}\to M^{p,q}} \, .
\end{equation}
  
We take some liberty in the interpretation of the operator norm
involved. As we have seen, the spaces we can chose depend mostly on
the quality of the window. 

Here are some precise approximation theorems. 

\begin{tm} \label{app}
Assume that $\cG (g, \Lambda )$ is a   (Parseval) Gabor frame for $\lrd $. 

  (i) If $g\in M^1(\rd )$ and  $\sigma \in \mif  $,
  then $E_N(\sigma ) \to 0$ 
 (in the operator norm on $M^{p,q}$ for all $1\leq p,q\leq \infty $).

(ii) If $g\in M^1_{v}$ and $\sigma \in
M^{\infty,\infty} _{v\circ j\inv}$, then  
$$
 E_N(\sigma )  \leq C \, \|\sigma \|_{M^{\infty,\infty} _{v\circ j\inv}} \, \sum
 _{|\nu | > N} v(\nu )\inv \, .
$$
  
(iii) If $g\in M^1_{v}$ and $\sigma \in
M^{\infty , 1} _{v\circ j\inv}$, then  
$$
 E_N(\sigma )   \leq C \,  \|\sigma \|_{M^{\infty ,1}
  _{v\circ j\inv}} \, \sup _{|\nu | >N} v(\nu )\inv \, .
$$

(iv) If $g\in M^1_{\langle \zeta \rangle ^s}$ and $\sigma \in
M^{\infty,\infty} _{\langle \zeta \rangle ^s}$ and 
$s>2d$, then  
$$
 E_N(\sigma )  \leq C \, \|\sigma \|_{M^{\infty,\infty} _{\langle
     \zeta \rangle ^s}} \, N^{2d-s}\, .
$$
  
(v) If $g\in M^1_{\langle \zeta \rangle ^s}$ and $\sigma \in
M^{\infty , 1} _{\langle \zeta \rangle ^s}$, then  
$$
 E_N(\sigma )   \leq C \, \|\sigma \|_{M^{\infty ,1}
  _{\langle \zeta \rangle ^s}} \, N^{-s} \, .
$$
\end{tm}

\begin{proof}
 (i)  By Proposition~\ref{triv} we represent $\sdd $ as a sum of Gabor multipliers 
$$
\sdd = \sum _{\nu \in \Lambda } \pi (\nu ) M_{\mathbf{a}_\nu } \, ,
$$
such that $\sum _{\nu \in \Lambda } \|\bba _\nu \|_\infty <
\infty$. Since the sum  converges in the operator norm on $M^{p,q}$,
the difference between $\sdd $ and its approximation  becomes  
$$
\sdd - \sum _{\nu \in \Lambda } \pi (\nu ) M_{\mathbf{a}_\nu } = \sum
_{|\nu |>N } \pi (\nu ) M_{\mathbf{a}_\nu } \, .
$$
 Taking operator norms,  we obtain the following
estimate for the error $E_N(\sigma )$:
\begin{eqnarray}
E_N(\sigma ) &=& \|\sdd - \sum _{\nu \in \Lambda } \pi (\nu )
M_{\mathbf{a}_\nu }\|_{M^{p,q} \to M^{p,q}} \notag \\
&=& \|\sum
_{|\nu |>N } \pi (\nu ) M_{\mathbf{a}_\nu } \|_{M^{p,q} \to M^{p,q}}
\notag \\
&\leq & \sum _{|\nu |>N } \|M_{\mathbf{a}_\nu } \|_{M^{p,q} \to
  M^{p,q}} \notag \\
&\leq & \sum _{|\nu |>N } \|\mathbf{a}_\nu \|_\infty \, . \label{eq:12}
\end{eqnarray}

 Since  $\sigma \in \mif $ and  $\sum _{\nu \in \Lambda }
 \|\mathbf{a}_\nu \|_\infty < \infty $,  \eqref{eq:12} implies that  $E_N(\sigma ) \to 0$. 

(ii) By Proposition~\ref{triv}  the assumption $\sigma \in M^{\infty ,\infty
}_{v\circ j\inv}  $ implies that $\|\mathbf{a}_\nu
\|_\infty \leq C \|\sigma \|_{M^{\infty, \infty } _{v\circ j\inv}}\,  v(\nu )\inv$. Consequently, by \eqref{eq:12} the
approximation error is at most
$$
E_N(\sigma ) \leq \sum _{|\nu |>N } \|\mathbf{a}_\nu \|_\infty \leq C
\|\sigma \|_{M^{\infty , \infty } _{v\circ j\inv}} \, \sum _{|\nu |>N }
  v(\nu )\inv  \, .
$$

If $v(z) = (1+|z|)^s$, then $\sum _{|\nu |>N} (1+|\nu|)^{-s} \leq C
N^{2d-s}$, whence assertion (iv).

(iii) Likewise, if $\sigma \in
M^{\infty , 1} _{v\circ j\inv}$, then by Theorem~\ref{sparsity} 
$$ \sum _{\nu \in \Lambda } \|\bba _\nu \|_\infty v(\nu )  \leq C
\|\sigma \|_{\mif _{v\circ j\inv} } < \infty  
$$
and thus 
$$
E_N(\sigma ) \leq \sum _{|\nu |>N } \|\mathbf{a}_n \|_\infty  \leq 
\sup _{|\nu |>N} v(\nu )\inv  \, \sum _{|\nu |>N }\|\bba \|_\infty \, v(\nu )\, .
$$

If $v(z) = (1+|z|)^s$, then $\sup _{|\nu |>N} v(\nu )\inv
\leq (1+N)^{-s}$, whence (v). 
\end{proof}

\rems\ 1. The error estimates of Theorem~\ref{app} yields an estimate
for the size of the cut-off parameter $N$ and thus of the number of
shifted \gmu s required for a good approximation. For instance, if
$\sigma \in \mif _{\langle \zeta \rangle ^s } (\rdd )$ and a tolerance $\epsilon
>0$ are given, then the error estimate
$$
 E_N(\sigma )   \leq C \, \|\sigma \|_{M^{\infty ,1}
  _{\langle \zeta \rangle ^s}} \, N^{-s} \, $$
implies that 
$$
N > \Big(\frac{(C\|\sigma \|_{\mif _{\langle \zeta \rangle ^s }}}{\epsilon} \Big)^{1/s}
\, .
$$
It is possible to say more about the dependency of  the constant $C$
 on   $g$ and  $\Lambda $  by precise bookkeeping  in the results
 about almost diagonalization of \psdo s in \cite{Gro06}.

2. One may wonder whether a converse of Theorem~\ref{app} holds
and whether the quality of approximation as expressed by
Theorem~\ref{app} characterizes the symbol classes. This guess is
false, because the operator norm $\|\sdd \|_{M^{p,q} \to M^{p,q}}$ and
$\|\sigma \|_{\mif} $ are not equivalent. Indeed, following an idea of 
Klotz~\cite{klotz} about approximation algebras of matrices,  one can
define a new  symbol class  directly  by  the 
approximation properties of \gmu s as follows: We say that $\sigma \in
\cA ^s(\rdd )$, if 
$$
\inf _{\bba _\nu \in \ell ^\infty (\Lambda )} \|\sdd - \sum _{\nu \in
  \Lambda , |\nu | >N} M_{\bba _\nu } \|_{L^2 \to L^2} \leq C N^{-s}
\, .
$$
  This class of symbols contains $M^{\infty ,\infty } _{\langle \zeta
    \rangle ^s}$. By comparison with the corresponding matrix algebras
  one can verify that $M^{\infty ,\infty } _{\langle \zeta
    \rangle ^s}$ is strictly smaller than  $\cA ^s(\rdd )$.

\section{Application}

Finally we mention the problem that has motivated the approximation
theorems in the previous section. This problem concerns the
transmission and decoding of  digital information by a variant of
(orthogonal) frequency division multiplexing (OFDM).

Here is a very coarse description of this procedure. 

1. Given is  ``digital information'' in the form of a finite sequence
$\{c_\lambda :  \lambda \in \Lambda \}$, where the data are taken from a
finite alphabet, usually just $c_\lambda \in \{-1,1\}$ or $c_\lambda
\in \{ \pm 1\pm i\}$.

2.  \emph{D/A conversion\/}:  To transmit these data, they are converted to an
analog signal  of the form 
$$
f= \sum _{\mu \in \Lambda } c_\mu 
\pi (\mu )g \, .
$$
for a suitable  pulse $g$.  Usually the lattice is taken to be a
rectangular lattice  $\Lambda =   \alpha \bZ \times \beta \bZ $. 
 Clearly, the coefficients must be uniquely determined by
 $f$. Therefore, and  in contrast to the previous discussions, the
 common assumption in wireless communication is that $\cG (g,\Lambda
 )$ is a Riesz sequence (for a closed subspace) in $L^2 (\bR )$.

3. \emph{Transmission of the signal $f$\/}: The  analog signal $f$ is
transmitted by a sender. This is a
physical process subject to  the laws of physics, in particular  the wave equation.

4. \emph{Distortion of $f$\/}: During the transmission the signal is distorted
by various effects. The most common effects are  time lags due to
reflection at obstacles and the Doppler effect due to relative motion
between transmitter and receiver. Thus  the distortion can be modelled
by  a superposition of \tfs s, and the   received signal is of the
form
\begin{equation}
  \label{eq:cd7}
  \tilde{f}(t) = \int _{\bR } \hat{\sigma} (\eta , u) M_\eta T_{-u}f(t)
\, du d\eta \, .
\end{equation}
Here the weighting factor $\hat{\sigma}$ models the  physical
details of the transmission, such as the reflectivities of obstacles
and the Doppler effect.
It is well known that the distortion~\eqref{eq:cd7} is precisely  the
\psdo\ with symbol 
$\sigma $, thus
\begin{equation}
  \label{eq:10}
\tilde{f} = \sdd f =  \sum _{\mu \in \Lambda } c_\mu \sdd (\pi (\mu
)g) \, .  
\end{equation}
Under the natural assumption of a  maximum Doppler shift $\nu _0$  and a maximum
time lag $\tau _0$, the weight function $\hat{\sigma }$ possesses a
compact support in $[0,\tau_0]\times [-\nu _0, \nu 
_0]$. If $\hat{\sigma }$ were  a bounded function, then  the
  distortion operator $\sdd $  would be a Hilbert-Schmidt operator and
 could not be invertible.  Therefore  Strohmer~\cite{str06} proposed
 the generalized Sj\"ostrand class 
  $M^{\infty , 1} _{e^{a|\zeta |}}$  with an exponential weight
  $v(\zeta ) = e^{a|\zeta|}$  as a 
  suitable symbol class  to model the distortion of time-varying systems in
  wireless communications. This class contains the distortion free
  channel corresponding to the identity operator (with symbol $\sigma
  \equiv 1$) and time-invariant channels corresponding to convolution
  operators.

5. \emph{ A/D conversion\/}: At the receiver we  decode the
original digital information $c_\lambda $ from the distorted
signal by   taking correlations with \tfs s as follows: 
\begin{eqnarray*}
y_\lambda &=&  \langle \tilde{f} , \pi (\lambda )g\rangle = \sum _{\mu
  \in \Lambda } c_\mu \underbrace{\langle \sdd (\pi (\mu )g), \pi (\lambda )
g\rangle} = (M(\sigma ) \mathbf{c})_\lambda  \, . 
  \end{eqnarray*}
The output vector  $\mathbf{y}= (y_\lambda )_{\lambda \in \Lambda }  $
is completely determined by the received signal 
$\tilde{f}$. With definition~\eqref{eq:cd8} the input-output relation can be
written as the infinite system of equations 
\begin{equation}
  \label{eq:11}
\mathbf{y} = M(\sigma )\mathbf{c}  
\end{equation}
In wireless communications the  matrix 
$M(\sigma ) $ is called the \emph{channel matrix}. 

For \eqref{eq:11} to be well posed,  the Gabor system $\cG (g, \Lambda )$ is assumed to
be a Riesz basis for its span $\cH _{g,\Lambda } = \mathrm{span} \, \cG
(g, \Lambda  )$,  and the distortion $\sdd $ is assumed to be an
invertible operator on $\lrd $. 
  Then the channel matrix $M(\sigma )$ is (boundedly) invertible on
  $\ell ^2(\Lambda )$, and the solution to \eqref{eq:11} is well
  defined. 
 
6. \emph{Equalization\/}: To solve for the original information $\mathbf{c}$ we
need to solve the system of equations~\eqref{eq:11} and find the
solution  $\mathbf{c} = M(\sigma )\inv \mathbf{y}$. 

 At this point occurs  amazing assumption that is  taken for granted by
engineers in this field: \emph{the channel matrix $M(\sigma )$ is
  assumed to be a diagonal matrix.}

With this assumption, the  distortion $\tilde{f}$ is simply 
\begin{eqnarray*}
 \tilde{f} = \sdd f &=& \sum _{\lambda \in \Lambda  } \langle \sdd \pi
 (\lambda ) g, \pi (\lambda )g\rangle \, \langle f, \pi  (\lambda
 )g\rangle \pi ( \lambda ) g \\
&=& \sum _{\lambda \in \Lambda  } M(\sigma ) _{\lambda \lambda } \, \langle f, \pi  (\lambda
 )g\rangle \pi ( \lambda ) g \, . 
\end{eqnarray*}
In other words, wireless communications  works with  the
implicit assumption that the  distortion operator    is a Gabor
multiplier with respect to a Riesz sequence $\cG (g,\Lambda )$. The
main motivation  for this 
assumption seems to be its  convenience and simplicity. With  this
assumption  the solution to~\eqref{eq:11} is
simply
$$
c_\lambda = M(\sigma )_{\lambda ,\lambda } \inv y_\lambda \, .
$$
One can show that the channel matrix with respect to a Gabor frame can
never be a diagonal matrix.  
However, if $\sigma \in
\mif _v$ with exponential weight $v$, then the channel matrix is
almost diagonal with exponential decay.  Theorem~\ref{app} then
guarantees  that the distortion operator $\sdd $ is 
approximated extremely well  by a finite number of shifted Gabor multipliers, or
equivalently, $M(\sigma )$ is approximated extremely well by a banded
matrix with few bands. In this case one may hope to improve the
accuracy of inversion of \eqref{eq:11} by using a banded approximation
of $M(\sigma )$ instead of the main diagonal only. In other words, we 
use an approximation of $\sdd $ by a finite sum of shifted Gabor
multipliers. 

This idea has been tested and implemented in
~\cite{HGM07,HGM09}. Combined with other tricks, 
the approximation of the distortion operator by a sum of \gmu s
contributed significantly  to a low complexity equalizer for
time-varying systems. This equalizer  performs best  in strongly
time-varying environments. 


\def\cprime{$'$} \def\cprime{$'$} \def\cprime{$'$} \def\cprime{$'$}
  \def\cprime{$'$}


\begin{thebibliography}{10}

\bibitem{bekka04}
B.~Bekka.
\newblock Square integrable representations, von {N}eumann algebras and an
  application to {G}abor analysis.
\newblock {\em J. Fourier Anal. Appl.}, 10(4):325--349, 2004.

\bibitem{BP06}
J.~J. Benedetto and G.~E. Pfander.
\newblock Frame expansions for {G}abor multipliers.
\newblock {\em Appl. Comput. Harmon. Anal.}, 20(1):26--40, 2006.

\bibitem{butzer}
P.~L. Butzer and R.~J. Nessel.
\newblock {\em Fourier analysis and approximation}.
\newblock Academic Press, New York, 1971.
\newblock Volume 1: One-dimensional theory, Pure and Applied Mathematics, Vol.
  40.

\bibitem{CG03}
E.~Cordero and K.~Gr{\"o}chenig.
\newblock Time-frequency analysis of localization operators.
\newblock {\em J. Funct. Anal.}, 205(1):107--131, 2003.

\bibitem{DT09}
M.~{D}{\"o}rfler and B.~{T}orresani.
\newblock {R}epresentation of operators in the time-frequency domain and
  generalized {G}abor multipliers.
\newblock {\em {J}. {F}ourier {A}nal. {A}ppl.},  2009, to appear.

\bibitem{feiSTSIP}
H.~G. Feichtinger.
\newblock Modulation spaces: looking back and ahead.
\newblock {\em Sampl. Theory Signal Image Process.}, 5(2):109--140, 2006.

\bibitem{fg89jfa}
H.~G. Feichtinger and K.~Gr{\"o}chenig.
\newblock Banach spaces related to integrable group representations and their
  atomic decompositions. {I}.
\newblock {\em J. Functional Anal.}, 86(2):307--340, 1989.

\bibitem{fg97jfa}
H.~G. Feichtinger and K.~Gr{\"o}chenig.
\newblock Gabor frames and time-frequency analysis of distributions.
\newblock {\em J. Functional Anal.}, 146(2):464--495, 1997.

\bibitem{FHK04}
H.~G. {F}eichtinger, M.~{H}ampejs, and G.~{K}racher.
\newblock {A}pproximation of matrices by {G}abor multipliers.
\newblock {\em {I}{E}{E}{E} {S}ignal {P}roc. {L}etters}, 11(11):883-- 886,
  {N}ovember 2004.

\bibitem{FN03}
H.~G. Feichtinger and K.~Nowak.
\newblock A first survey of {G}abor multipliers.
\newblock In {\em Advances in Gabor analysis}, Appl. Numer. Harmon. Anal.,
  pages 99--128. Birkh\"auser Boston, Boston, MA, 2003.

\bibitem{book}
K.~Gr{\"o}chenig.
\newblock {\em Foundations of time-frequency analysis}.
\newblock Birkh\"auser Boston Inc., Boston, MA, 2001.

\bibitem{Gro06ped}
K.~Gr{\"o}chenig.
\newblock A pedestrian's approach to pseudodifferential operators.
\newblock In {\em Harmonic analysis and applications}, Appl. Numer. Harmon.
  Anal., pages 139--169. Birkh\"auser Boston, Boston, MA, 2006.

\bibitem{Gro06}
K.~Gr\"ochenig.
\newblock Time-frequency analysis of {S}j\"ostrand's class.
\newblock {\em Revista Mat. Iberoam.}, 22(2):703--724, 2006.


\bibitem{GR08}
K.~Gr{\"o}chenig and Z.~Rzeszotnik.
\newblock Banach algebras of pseudodifferential operators and their almost
  diagonalization.
\newblock {\em Ann. Inst. Fourier (Grenoble)}, 58(7):2279--2314, 2008.

\bibitem{HGM09}
M.~{H}ampejs, P.~{S}vac, G.~{T}aub{\"o}ck, K.~{G}r{\"o}chenig, F.~{H}lawatsch,
  and G.~{M}atz.
\newblock Sequential {L}{S}{Q}{R}-based {I}{C}{I} equalization and decision
  feedback {I}{S}{I} cancelalation in pulse-shaped multicarrier systems.
\newblock {P}roc. {I}{E}{E}{E} {S}{P}{A}{W}{C}09, Helsinki.

\bibitem{hormander3}
L.~H{\"o}rmander.
\newblock {\em The analysis of linear partial differential operators. {III}},
  volume 274 of {\em Grundlehren der Mathematischen Wissenschaften [Fundamental
  Principles of Mathematical Sciences]}.
\newblock Springer-Verlag, Berlin, 1994.
\newblock Pseudo-differential operators, Corrected reprint of the 1985
  original.

\bibitem{klotz}
A.~Klotz.
\newblock {\em Noncommutative approximation: smoothness and approximation and
  invertibility in {B}anach algebras}.
\newblock PhD thesis, Univ. Vienna, 2009.

\bibitem{Sjo94}
J.~Sj{\"o}strand.
\newblock An algebra of pseudodifferential operators.
\newblock {\em Math. Res. Lett.}, 1(2):185--192, 1994.

\bibitem{Sjo95}
J.~Sj{\"o}strand.
\newblock Wiener type algebras of pseudodifferential operators.
\newblock In {\em S\'eminaire sur les \'Equations aux D\'eriv\'ees Partielles,
  1994--1995}, pages Exp.\ No.\ IV, 21. \'Ecole Polytech., Palaiseau, 1995.

\bibitem{str06}
T.~Strohmer.
\newblock Pseudodifferential operators and {B}anach algebras in mobile
  communications.
\newblock {\em Appl. Comput. Harmon. Anal.}, 20(2):237--249, 2006.

\bibitem{HGM07}
G.~{T}aub{\"o}ck, M.~{H}ampejs, G.~{M}atz, F.~{H}lawatsch, and
  K.~{G}r{\"o}chenig.
\newblock {L}{S}{Q}{R}-based {I}{C}{I} equalization for multicarrier
  communications in strongly dispersive and highly mobile environments.
\newblock {P}roc. {I}{E}{E}{E} {S}{P}{A}{W}{C}07.

\bibitem{Toft02a}
J.~Toft.
\newblock Continuity properties in non-commutative convolution algebras, with
  applications in pseudo-differential calculus.
\newblock {\em Bull. Sci. Math.}, 126(2):115--142, 2002.

\bibitem{wong02}
M.~W. Wong.
\newblock {\em Wavelets Transforms and Localization Operators}, volume 136 of
  {\em Operator Theory Advances and Applications}.
\newblock Birkhauser, 2002.

\end{thebibliography}

\end{document}